\documentclass[11pt]{article}
\usepackage{amssymb}
\usepackage{latexsym}
\newcommand{\reals}{{\mathbb R}}
\newcommand{\torus}{{\mathbb T}}

\newcommand{\calg}{{\cal G}}

\newcommand{\calt}{{\cal T}}

\newcommand{\calu}{{\cal U}}
\newcommand{\del}{\partial}

\newcommand{\half}{\textstyle{\frac{1}{2}}}
\newcommand{\frakg}{\mathfrak{g}}
\newcommand{\frakh}{\mathfrak{h}}
\newcommand{\frakm}{\mathfrak{m}}
\newcommand{\arrows}{\,\lower1pt\hbox{$\longrightarrow$}\hskip-.24in\raise2pt
             \hbox{$\longrightarrow$}\,}
\title{{\bf Linearization Problems for Lie Algebroids and Lie Groupoids
}\\In memory of Mosh\'e Flato, 1938-1999} \author{Alan
Weinstein\thanks{Research partially supported by NSF Grants
DMS-96-25122 and DMS-99-71505
and the Miller Institute for Basic Research in Science.}
\\Department of Mathematics\\ University of California\\ Berkeley, CA
94720 USA\\ {\small(alanw@math.berkeley.edu)}}
\date{November 1999}
\begin{document}
\setlength{\baselineskip}{15pt}

\maketitle

\section{Introduction}
\label{sec-intro}
Why is it so
hard to prove the linearizability of Poisson structures with
semisimple linear part?  Conn published proofs 
about 15 years ago in a pair of papers \cite{co:analytic}\cite{co:smooth}  full of
elaborate estimates.  Except for the somewhat more conceptual
reformulation by Desolneux-Moulis \cite{de:linearisation} in the
smooth case, no 
simplification of Conn's proofs has appeared.   

This is a mystery to me, because analogous theorems about the
linearizability of actions of semisimple groups near their fixed
points were proven in the compact (smooth or analytic) case
by Bochner \cite{bo:compact} using a
simple averaging method  and
in the noncompact analytic case by Guillemin--Sternberg
\cite{gu-st:remarks} and Kushnirenko \cite{ku:analytic}, who used analytic
continuation from 
the compact case--a nonlinear version
of ``Weyl's unitary trick''.  Hermann \cite{he:formal} 
established formal linearization for actions of general semisimple
algebras, using cohomological methods similar to those which will
appear several times in the present report.

After Conn's work appeared, I tried without success to prove his results by
simple averaging.  In a conversation
in 1994 over coffee
in the Jardin du Luxembourg, Mosh\'e Flato and Daniel Sternheimer
revived my interest in the subject. Their
experience with Jacques Simon on other linearization problems
(see for instance \cite{fl-si:linearization})
led them to expect
simpler proofs of Conn's theorems.  In
addition, they pointed out that Poisson linearization
could be seen as an infinite-dimensional Levi decomposition. 
Although I
expressed strong skepticism at the time, Mosh\'e and Daniel's
optimism has remained in my mind and kept me from abandoning the
problem.  In particular, it has resulted in new investigations which,
although they have not produced a Poisson linearization theorem (and
maybe never will!), have led to a few results and many questions which
are interesting in their own right.  

These new results and problems form the content of the present paper,
which I offer as a small memorial to Mosh\'e.  The paper also serves
to record a talk which I gave at the Dijon Mathematical Physics
Seminar in June, 1999, at Daniel's invitation.\footnote{An earlier
version of this paper was presented in a lecture at the Chern
Symposium at MSRI \cite{we:from}.}  The occasion was much saddened by
Mosh\'e's absence, but it has been gratifying to see the Seminar
continue as an ongoing international ``workshop'' for discussion of
many of Mosh\'e's favorite scientific issues.

\section{Linearization and Levi decomposition}
\label{sec-levi}
The linearization of a Poisson structure with semisimple linear part
can be seen, as mentioned
above, as something like a Levi decomposition of a Lie algebra of functions
with the Poisson bracket operation.  Here are the details.

Let $\pi=\half\sum(\sum c_{ij}^k x_k +{\bf O}(x^2))\frac{\del}{\del
x_i}\wedge\frac{\del}{\del x_j} $ be a Poisson structure defined on a
neighborhood $\calu$ of the origin in $\reals^n$, i.e. $\{x_i,x_j\}=
\sum c_{ij}^k x_k + {\bf O}(x^2).$ In the Lie algebra of germs at $0$
of smooth functions, the germs of functions vanishing at $0$ form a
Lie subalgebra $\frakm$ in which those vanishing to order at least 2
form a Lie algebra ideal $\frakm ^2$.  The quotient $\frakm/\frakm^2 =
\frakg$ may be identified with the cotangent space of $\reals^n$ at
$0$, the (finite dimensional) {\bf cotangent Lie algebra}
 whose structure constants are
just the first Taylor coefficients $c_{ij}^k$ of the Poisson structure
$\pi$.  The {\bf linearization problem} is to find (perhaps after
shrinking the neighborhood $\calu$) new coordinates $(y_1,\ldots,y_n)$
centered at $0$ such that $\{y_i,y_j\}= \sum c_{ij}^k y_k$ without any
higher order terms.

To view a solution of the linearization problem for semisimple
$\frakg$ as a Levi decomposition, we note first that, in the exact
sequence of Lie algebras
$$0\rightarrow \frakm ^2\rightarrow \frakm\rightarrow \frakg\rightarrow 0,$$
the quotient $\frakg$ is (by assumption) semisimple, while the kernel
$\frakm^2$ has a certain nilpotency property.  
By this last statement, we mean that
the spaces $\frakm^k$ of functions vanishing to order at least $k$
have the property that $[\frakm^2,\frakm^k] \subseteq \frakm^{k+1}$.
Unfortunately, the intersection 
$\frakm^{\infty}=\cap_{k}\frakm^k$ is not zero, but consists of the
germs of functions which vanish to infinite order at the
origin.\footnote{There are some interesting questions to be answered
about $\frakm^{\infty}$.  Does it have any semisimple quotients?  Does
it have any finite-dimensional quotients at all?   Equivalently, does
it have any finite-dimensional representations?  What are the its
finite-dimensional subalgebras?}

A ``Levi decomposition'' of $\frakm$ should be a vector space direct sum
decomposition $\frakm = \frakh \oplus \frakm^2$ in which $\frakh$ is a
subalgebra or, equivalently, 
a Lie algebra homomorphism $\phi:\frakg\to \frakm$ which
splits the exact sequence.

Given a linearization $(y_1,\ldots,y_n)$ of the Poisson structure, the
subalgebra spanned by the germs of the $y_i$ is the required
$\frakh$.  Conversely, given a Levi decomposition and letting
$(\eta_1,\ldots,\eta_n)$ be the preimages in $\frakh$ of the basis
$(x_1+\frakm^2,\ldots,x_n + \frakm^2)$ of $\frakg$, representatives $y_i$ of the
germs $\eta_i$ have a common domain on which they solve the linearization
problem.  

Note that everything above remains true if the
algebra of smooth germs is replaced by either that of real analytic germs or
that of formal power series, with the additional feature that $\frakm^{\infty}$
is now reduced to zero.

Unfortunately, I do not know of any proof of the standard Levi decomposition
theorem which would apply here, 
so it seems that we have only  reformulated the linearization
problem.  Nevertheless, this reformulation leads to some interesting new
questions.

Finally, we note that A.Wade \cite{wa:normalisation} has
found a formal normalization of Poisson structures at general singular
points which is based on the Levi decomposition of the cotangent Lie
algebra itself.  

\section{Lie algebroids}
\label{sec-algebroids}

The linearization problem for Poisson manifolds has analogues for
Lie algebroids and Lie groupoids.  To describe them, we recall
that a Lie algebroid over a manifold $M$ is a vector bundle
$A\rightarrow M$ with  a Lie
algebra structure (over $\reals$) on its space $\Gamma (A)$ of
smooth sections  and a bundle map $\rho:A\rightarrow TM$
inducing a Lie algebra homomorphism from sections of $A$ to
vector fields on $M$, such that $[a,fb]=f[a,b]+(\rho(a)\cdot f)b$
for sections $a$ and $b$ and functions $f$.  Locally (which is all we
need for this paper), a Lie algebroid over an
open subset $\calu$ of $\reals^n$ is specified by the bracket
relations $[e_i,e_j]=\sum c_{ij}^k (x)e_k$ among  a basis of sections
of $A$ and the component description $\rho(e_i)=\sum
b_{ij}(x)\frac{\del}{\del x_j}$ of the anchor map.  The axioms of a
Lie algebroid then become differential equations relating the
{\bf structure functions} $c_{ij}^k$ and $b_{ij}$.  

As with Poisson manifolds, the local classification of Lie algebroids
can be reduced by a splitting theorem (modeled on the splitting
theorem of Dazord \cite{da:feuilletages} for singular foliations) to the
totally singular case where the anchor map $\rho$ vanishes at a
point.  In coordinates centered at such a point, we may write
$c_{ij}^k (x) = c_{ij}^k(0) + {\bf O}(x)$ and
$b_{ij}(x)=\frac{\del b_{ij}}{\del x_k}(0) x_k + {\bf O}(x^2).$
Our problem is to  present the Lie algebroid
in such a way that the higher order 
terms in the structure functions disappear.  Although we have more
conditions to satisfy than in the Poisson linearization problem, we
also have more variables at our disposal, since we can change the basis
of sections as well as the coordinates on the base.  

We first  try to change the basis  to reduce the
functions $c_{ij}^k (x)$ to the constants $c_{ij}^k (0).$  This is
equivalent to making the Lie algebroid into one of a special kind.
 Namely, the constants $c_{ij}^k (0)$
define a Lie algebra structure on the fibre of $A$ at $0$, and the
vector fields $\rho(e_i)$ give an action of this Lie algebra on
$\calu$.  Conversely,  when a Lie algebra $\frakg$ acts on a manifold
$M$, there is a Lie algebroid structure on the trivial bundle
$\frakg\times M$ for which the anchor is given by the action and the
constant sections form a subalgebra on which the bracket is that of $\frakg$.  So our
task is to determine when a given Lie algebroid is, near a point where
the anchor vanishes, such an {\bf action Lie algebroid}.  If this task is
accomplished, we can finish simplifying the Lie algebroid structure
functions by linearizing the action, under the usual hypotheses of
semisimplicity and, in the smooth setting, compact type.  

Making $A$ into an action Lie algebroid on a neighborhood amounts
precisely to finding a Lie algebra homomorphism which splits the exact
sequence:
$$0\rightarrow\mbox{sections vanishing at $0$}\rightarrow
\mbox{sections of $A$}\rightarrow \mbox{fibre at $0$}\rightarrow 0,$$
which is analogous to the sequence in the Poisson linearization problem.
In particular, if we pass to germs, the kernel of the restriction map
is once again ``topologically nilpotent'' in the sense discussed
in Section \ref{sec-levi}.  The analogy is so close that it is
tempting to transfer the proofs of Poisson linearization theorems to
the Lie algebroid case.  In the formal category, this is no problem.
The obstructions to stepwise lifting from the fibre $\frakg$ at $0$ to higher
and higher jets of sections lie in the spaces
$H^2(\frakg,S^k(\reals^n)\otimes \frakg)$, where $\frakg$ acts on
itself by the adjoint representation and on the symmetric tensor power
$S^k(\reals^n)$ via the action on $\reals^n$ given by
$e_i\mapsto\frac{\del b_{ij}}{\del x_k}(0).$
These cohomology spaces vanish when $\frakg$ is
semisimple.\footnote{For formal Poisson linearization, the relevant
cohomology spaces are $H^2(\frakg,S^k(\frakg)).$}

It seems likely that Conn's methods can be extended to establish
linearizability in the smooth and analytic settings.  I have not yet tried
to do this, however, preferring another approach using
averaging which might eventually lead to a new proof of Poisson linearization.
The following sections describe this approach.

\section{Lie groupoids}
\label{sec-groupoids}
Averaging requires an action of a compact Lie group, not just of its
Lie algebra.  We therefore pose a linearization problem for Lie
groupoids, which are the ``integrated'' form of Lie algebroids (see
\cite{ca-we:geometric} or \cite{ma:lie}).  Note, though, that not every
Lie algebroid is integrable to a groupoid \cite{al-mo:suites}, so that
linearization of Lie groupoids will linearize only certain Lie
algebroids.

Let $\Gamma$ be a Lie groupoid over $M$
with source and target maps $\beta$ and $\alpha$; the product $gh$ is
defined when $\beta(g)=\alpha(h)$.  $\Gamma$ is a {\bf proper
groupoid} if $(\alpha,\beta):\Gamma\rightarrow M\times M$ is a proper
mapping and if the source and target maps are locally trivial
fibrations.  This term, like many other terms in groupoid theory,
comes from the following example.

If $G$ is a group acting on $M$, then $G\times M$ is a Lie groupoid
over $M$ with $\alpha(g,x)=gx$, $\beta(g,x)=x$, and product
$(g,hx)(h,x)=(gh,x)$.  We call this the {\bf action groupoid}
associated to the group action.  It is a proper groupoid if and only
if the action is a proper action.  

For a general groupoid $\Gamma$ over $M$, the {\bf isotropy group}
$\Gamma_x$ of $x$ in $M$ is defined as
$\alpha^{-1}(x)\cap\beta^{-1}(x).$  Since $\Gamma_x =
(\alpha,\beta)^{-1}(x,x)$, every isotropy group in a proper groupoid is
compact.  

A subset $\calu \subseteq M$ is called {\bf invariant} under $\Gamma$
if $\alpha(\beta^{-1}(u))\subseteq \calu$.  An invariant point is
called a {\bf fixed point}.

We show in \cite{we:proper} that every neighborhood of a fixed point of a
proper groupoid contains an invariant neighborhood, so that the study of proper
groupoids can be localized around fixed points.  We make the following
conjecture.

\bigskip
\noindent
{\bf Proper Groupoid Structure Conjecture}.  {\em  If $x$ is a fixed
point of a proper groupoid $\Gamma$, then $x$ has an invariant
neighborhood on which $\Gamma$ is isomorphic to the action groupoid
associated to some action of the compact group $\Gamma_x$.}
\bigskip

If the ``PGS'' conjecture above is true, the next step is to linearize
the $\Gamma_x$ action, which would give an equivalence between proper
groupoids near their fixed points and linear actions of compact
groups.  From here it should be possible, using a slice theorem, to
get a normal form for a proper groupoid in the neighborhood of any
orbit $\beta (\alpha^{-1}(x))$.  This program is described in \cite{we:proper}.

The PGS conjecture fits into the Levi
decomposition picture presented in Section \ref{sec-levi} via the
 exact sequence of groups:
$$1\rightarrow \calg_{\calu}^1\rightarrow
\calg_{\calu}\rightarrow\Gamma_x\rightarrow 1.$$  
Here, $\calg_{\calu}$ is the group of {\bf admissible sections} 
of $\Gamma$ over the invariant set $\calu$, i.e. the smooth
submanifolds $\sigma$ of the restricted groupoid $\Gamma_{\calu}$ for
which the restrictions $\alpha|_{\sigma}$ and $\beta|_{\sigma}$ are
diffeomorphisms from $\sigma$ to $\calu$.  The map
$\calg_{\calu}\rightarrow\Gamma_{x}$ is evaluation at $x$, and the
kernel $\calg_{\calu}^1$ consists of those sections which meet the
unit section at $x$.  Finding an isomorphism of $\Gamma_{\calu}$ with
an action groupoid amounts to finding a cross section
$\Gamma_x\rightarrow \calg_{\calu}$ which is a group homomorphism.  

Since the properness of $\Gamma$ implies that $\Gamma_x$ is compact,
we could construct a homomorphic section by averaging if the kernel
$\calg_{\calu}^1$ were the additive group of a vector space.  Of
course, this is not the case, but we do get a nice composition series
for $\calg_{\calu}^1$ if we pass to infinite jets at $x$ of admissible
sections of $\Gamma$.  We can then use the standard stepwise proof
(this time using group cohomology) to find an action groupoid
structure for $\Gamma$ over a ``formal neighborhood of $x$''.  Note
that semisimplicity plays no role here--compactness of $\Gamma_x$ is
enough, so it could be a torus, for instance.

In the analytic or smooth categories, we could try to imitate Conn's
proofs, but instead we propose another approach.  We confine our
discussion to the smooth case.  

Since $\alpha$ is
a locally trivial fibration, we can choose a cross
section $\sigma_0:\Gamma_x\rightarrow \calg_{\calu}$ for sufficiently
small $\calu$ such that the corresponding map from $\Gamma_x \times
\calu$ to $\alpha^{-1}(\calu)$ is a diffeomorphism.  This $\sigma_0$ will
not generally be a group homomorphism, so we will try to ``improve''
it to become one.
By the formal normal form result, we can assume (here we imitate Conn
\cite{co:smooth}) that $\sigma_0$ is a homomorphism modulo flat sections;
i.e. for each $g$ and $h$ in $\Gamma_x$, $\sigma_0(g)\sigma_0(h)$ and
$\sigma_0(gh)$ are tangent to infinite order at the point $gh$ where
they meet $\Gamma_x$.  Using our local trivialization of $\alpha$ and
local coordinates on a sufficiently small open neighborhood $\calu$ of
$x$, we can assume that, for all $g$ and $h$ in $\Gamma_x$,
$\sigma_0(g)\sigma_0(h)$ and $\sigma_0(gh)$ considered as maps from
$\calu$ to $\Gamma_{x}$ are as close as we wish, together with as many
derivatives as we desire.

Thus, giving $\calg_{\calu}$ a $C^k$-metric, we can say that
$\sigma_0$ is an ``almost homomorphism'' in the sense that
$d(\sigma_0(g)\sigma_0(h),\sigma_0(gh)$ is small for all $g$ and $h$
in $\Gamma_x$, and our problem is to approximate 
this almost homomorphism by  a homomorphism.

\section{Almost homomorphisms and almost invariant submanifolds}
\label{sec-almost}

There are quite a few theorems which assert that an almost
homomorphism is near a homomorphism (a very general formulation of
this problem was proposed by Ulam \cite{ul:sets}), but none of them suits our needs,
so we must
prove new ones.

The following two theorems are closest to what we are looking for.

\bigskip
\noindent
{\bf Almost Homomorphism Theorem.} (Grove--Karcher--Ruh
\cite{gr-ka-ru:group}).  {\em Let 
$G$ and $H$ be compact Lie groups.  Choose a bi-invariant metric on
$H$ such that: (1) the exponential map is an embedding when restricted
to the open ball of radius $\pi$ in the Lie algebra $\frakh$; (2) the
Lie bracket on $\frakh$ satisfies $||[v,w]|| \leq ||v||~||w||$.  (Such
a metric always exists.)  Let $\sigma_0:G\to H$ be a continuous map
such that $d(\sigma_0(g)\sigma_0(h),\sigma_0(gh))\leq q\leq \pi/6$ for
all $g$ and $h$ in $G$.  Then there is a (continuous) homomorphism
$\sigma:G\to H$ such that $d(\sigma_0(g),\sigma(g)) < 1.36 q$ for all
$g$ in $G$. 
}

\bigskip
\noindent
{\bf Almost Representation Theorem}.  (de la Harpe--Karoubi
\cite{de-ka:representations}) 
{\em Let $T_0$ be a continuous map from a compact group $G$ to the
group $H$ of invertible bounded linear operators on some Hilbert
space.  Let $K\geq 1$ and $\epsilon\leq 2^{-6}$ be real numbers such
that $||T_0(gh)-T_0(g)T_0(h)|| \leq  \epsilon (2K)^{-9}$ for all $g$
and $h$ in $G$.  Then there is a continuous homomorphism $T:G\to H$
such that $||T_0(g)-T(g)|| \leq \epsilon$ for all $g$ in $G$.
}

An intermediate result between those theorems and the one we want for
groupoid linearization would be the following, stated somewhat
imprecisely.

\bigskip
\noindent
{\bf Almost Action Conjecture}. {\em Let $G$ be a compact group and $H$
the group of diffeomorphisms of a compact manifold $M$.  If
$\sigma:G\rightarrow H$ is a map such that
the distance $d(\sigma_0(g)\sigma_0(h),\sigma_0(gh)$ is sufficiently 
small for all
$g$ and $h$
in $G$, then there is a homomorphism $\psi$ close to $\sigma$.  Here,
the distance between two diffeomorphisms is taken to be their $C^1$
distance defined using a riemannian metric $\mu$, and the measure of
``sufficient smallness'' will depend on the size of the first and
second derivatives of the maps in $\phi(G)$, as well as on upper bounds on
the absolute value of the curvature and the reciprocal of the
injectivity radius for the metric $\mu$. }

Even this conjecture is currently beyond our reach.  When
$G$ is finite, it should follow from the 
following result in \cite{we:almost}.  

\bigskip
\noindent
{\bf Almost Invariant Submanifold theorem}.  {\em Let $N$ be a compact
submanifold of a riemannian manifold $M$, and let $G$ be a compact
group acting on $M$.  Choose an invariant metric on $M$ such that: (1)
the exponential map of $M$, restricted to the normal bundle of $N$, is
(defined and) an embedding on the open ball bundle $B$ of radius 1;
(2) all the sectional curvatures of $M$ on $\exp(B)$ have absolute
value less than 1; (3) the exponential map of $M$ is (defined and) an
embedding on the ball of radius 1 in each tangent space of $\exp(B)$.
(This is always possible: multiply any invariant metric by a suitably
large constant.)  If the $C^1$ distance (defined below) $d(N,gN)$ is
less than $\epsilon < \frac{1}{20000}$ for each $g\in G$, then there
is a $G$-invariant submanifold $\overline{N}\subseteq M$ such that
$d(N,\overline{N})<136\sqrt{\epsilon}.$ 
}
\bigskip

The distance $d(N,N')$ between submanifolds in the theorem above is
defined as follows.  First of all, we assume that $N'$ is the image
under $\exp$ of a section $s$ of the normal bundle of $N$.  For each
$x\in N$, we take the maximum $a(x)$ of the following two numbers: the length
of the geodesic segment $\sigma$ from $x$ to $s(x)$, and the maximum
angle between unit vectors in $T_x N$ and $T_{x'}N'$, where vectors in
one space are moved to the other by parallel transport along
$\sigma$.  Now $d(N,N')$ is defined as $\max_{x\in N} a(x)$. 

Like   the proofs of  the    Almost  Homomorphism Theorem  and the  Almost
Representation   Theorem, the proof of  the Almost Invariant
Submanifold Theorem uses averaging
over the group  $G$.  The estimates involved  are nontrivial, but they
are somewhat simpler, and  certainly more geometric,  than the ones in
Conn's proof.  The challenge  now is to  work back from here to  prove
some linearization theorems without adding too much more complication.

\section{Properness and convexity}
Our interest in proper groupoids was originally motivated, not by
linearization problems, but by an attempt to understand the convexity
theorems of Atiyah \cite{at:convexity}, Guillemin--Sternberg
\cite{gu-st:convexity} \cite{gu-st:convexity2} and Kirwan
\cite{ki:convexity} as results in Poisson geometry.  Their theorems establish
convexity properties of the image of the momentum map $J:S\to
\frakg^*$ for a hamiltonian action of a compact Lie group $G$ on a
symplectic manifold $S$.  When $G$ is a torus, the theorem states that
$J(S)$ is a convex polyhedron, while for general compact $G$ it is
the intersection of $J(S)$ with a positive Weyl chamber which is
convex.  

When $G$ is simply connected, a momentum map is simply a Poisson map
from $S$ to $\frakg^*$, which suggests that there might be a convexity
theorem for Poisson maps $J:S\to P$ from symplectic manifolds to a
wider class of Poisson manifolds.  The case of the torus $\torus ^n$
shows that the problem involves more than just Poisson manifolds: in
this case, $P={\calt^n}^*$ has the zero Poisson structure, so that the
set of Poisson maps from $S$ to $P$ is closed under arbitrary
diffeomorphisms of $P$, which usually destroy any convexity properties
of the image. The point here is that $J$ must the momentum map for an
action of $\torus^n$ and not just of its universal covering
$\reals^n$.  

To express in geometric terms the choice of group associated to a
given (dual of a) Lie algebra, we recall that hamiltonian actions of
the Lie group $G$ correspond to symplectic actions of the symplectic
groupoid $T^*G\arrows \frakg^*$ \cite{mi-we:moments}. (The essentials of this
result were presented in \cite{we:lectures}, before the advent of symplectic
groupoids in the 1980's.)  We therefore pose the following conjecture
(having well in mind that some further hypotheses may be needed.)

\bigskip
\noindent
{\bf Poisson Convexity Conjecture.}  {\em Let $\Gamma\arrows P$ be a
proper symplectic groupoid over the Poisson manifold $P$, and let
$J:S\rightarrow P$ be the moment map (see \cite{mi-we:moments}) 
for an action of $\Gamma$
on a symplectic manifold $S$.  Then the image $J(S)\subseteq P$ has a
convexity property with respect to some affine structure attached to
the groupoid $\Gamma$.}  

\bigskip
The problem here is to establish the appropriate notion of convexity,
which requires an understanding of the structure of proper symplectic
groupoids.  Such an understanding would also tell us the extent to
which our conjecture goes beyond the known convexity theorems for
group actions.  

\bigskip
\noindent
{\bf Examples.} 
Let $G$ be a semisimple Lie group of noncompact type.  The {\bf
elliptic} subject $E \subset \frakg^*$ is defined to consist of those
elements whose coadjoint isotropy group in $G$ is compact.  (Under the
orbit method, $E$ corresponds to the discrete series of
representations of $G$.)  It appears
that $E$ is an open subset on which the restriction of
the coadjoint representation is a proper action.  By restricting the
symplectic groupoid $T^*G\arrows \frakg^*$ to $E$, we obtain a proper
symplectic groupoid for $E$.  Now a hamiltonian action of $G$ on $S$
may be called ``elliptic'' if the image $J(S)$ of the momentum map
lies in $E$.  A special case of our problem would be to establish a
convexity theorem for the momentum maps of such elliptic actions.

The previous example is still a transformation groupoid.  Another
example of a proper symplectic groupoid is the fundamental groupoid
$\pi(M)$ of
a compact, connected, symplectic manifold $M$ with finite fundamental group.
Symplectic actions of this groupoid correspond to symplectic actions
of the fundamental group.  (See Section 7.6 of
\cite{ca-we:geometric}.)  In this case, there is no convexity problem,
since momentum maps to the symplectic manifold $M$ are surjective, but
perhaps products $\pi(M) \times T^*G \arrows M \times \frakg^*$ will
give interesting examples.


\begin{thebibliography}{99}

\bibitem{al-mo:suites} 
Almeida, R., and Molino, P., Suites d'Atiyah et feuilletages
transversalement complets, {\em C.R. Acad. Sci. Paris} {\bf 300}
(1985), 13-15.

\bibitem{at:convexity}
Atiyah, M.F., Convexity and commuting Hamiltonians, {\em Bull. London
Math. Soc.} {\bf 14} (1982), 1-15.

\bibitem{bo:compact}
Bochner, S., Compact groups of differentiable transformations,
{\em Annals of Math.} {\bf 46} (1945), 372--381.

\bibitem{ca-we:geometric}
Cannas da Silva, A., and Weinstein, A., {\em 
Geometric Models for Noncommutative Algebras}, Berkeley
Math. Lecture Notes, Amer. Math. Soc., Providence, 1999.


\bibitem{co:analytic}
Conn, J.F., Normal forms for analytic Poisson structures, {\em Annals
of Math.} {\bf 119} (1984), 576-601.
 
\bibitem{co:smooth}
Conn, J.F., Normal forms for smooth Poisson structures, {\em Annals
of Math.} {\bf 121} (1985), 565-593.

\bibitem{da:feuilletages} 
Dazord, P., Feuilletages \`{a} singularit\'{e}s, {\em Nederl. Akad. Wetensch.
Indag. Math.} {\bf 47} (1985), 21-39.

\bibitem{de-ka:representations}
de la Harpe, P., and Karoubi, M., Repr\'esentations approch\'ees
d'un grope dans une alg\`ebre de Banach, {\em Manuscripta Math.}
{\bf 22} (1977), 293-310.

\bibitem{de:linearisation}
Desolneux-Moulis, N.,
Lin\'earisation de certaines structures de Poisson de classe
$C^{\infty}$. S\'eminaire de 
g\'eom\'etrie, 1985--1986, 
{\em Publ. D\'ep. Math. Nouvelle S\'er.} B, 86-4, 
Univ. Claude-Bernard, Lyon (1986), 55--68.

\bibitem{fl-si:linearization}
Flato, M.; Simon, J. 
On a linearization program of nonlinear field equations. 
Phys. Lett. B 94 (1980), no. 4, 518--522. 

\bibitem{gr-ka:how}
Grove, K., and Karcher, H.,
How to conjugate $C^1$-close group actions, {\em Math. Z.} {\bf 132}
(1973), 11-20.

\bibitem{gr-ka-ru:group}
Grove, K., Karcher, H., and Ruh, E.A. Group actions and
curvature, {\em Invent. Math.} {\bf 23} (1974), 31-48.

\bibitem{gu-st:remarks}
Guillemin, V. W., and  Sternberg, S.,
Remarks on a paper of Hermann,
{\em Trans. Amer. Math. Soc.} {\bf 130} (1968), 110--116. 


\bibitem{gu-st:convexity} 
Guillemin, V., and Sternberg, S., Convexity
properties of the moment mapping, {\em Invent. Math.} {\bf 67} (1982),
491-513.

\bibitem{gu-st:convexity2} Guillemin, V., and Sternberg, S., Convexity
properties of the moment mapping II, {\em Invent. Math.} {\bf 77} (1984),
533-546.

\bibitem{he:formal}
Hermann, R.,  The formal linearization of a semisimple Lie algebra
of vector fields
about a singular point, {\em Trans. Amer. Math. Soc.} {\bf 130} (1968) 105--109.

\bibitem{ki:convexity} 
Kirwan, F., Convexity
properties of the moment mapping III, {\em Invent. Math.} {\bf 77} (1984),
547-552.

\bibitem{ku:analytic}
Kushnirenko, A.G.,
Linear-equivalent action of a semisimple Lie group in the neighborhood
of a fixed point, {\em Funct. Anal. Appl.} {\bf 1} (1967), 89-90.

\bibitem{ma:lie} 
Mackenzie, K., {\em Lie Groupoids and Lie Algebroids in Differential
Geometry}, LMS Lecture Notes Series, {\bf 124}, Cambridge Univ. Press, 1987.

\bibitem{mi-we:moments} 
Mikami, K., and Weinstein, A., Moments and reduction for symplectic
groupoid actions, {\em Publ. RIMS Kyoto Univ.} {\bf 24}
(1988), 121--140.

\bibitem{ul:sets}
Ulam, S., {\em Sets, Numbers, and Universes}, MIT Press, Cambridge,
1974.

\bibitem{wa:normalisation}
Wade, A., Normalisation formelle de structures de Poisson,
{\em C. R. Acad. Sci. Paris S\'er. I Math.} {\bf 324} (1997), 531--536.

\bibitem{we:lectures} 
Weinstein, A., {\em Lectures on Symplectic Manifolds}, Regional
conference series in mathematics {\bf 29}, (American Mathematical Society,
Providence, 1977).

\bibitem{we:from}
Weinstein, A., From Riemann geometry to Poisson geometry and back
again, Lecture at Chern Symposium, MSRI, 1998, available at
{\tt http://msri.org/publications/video/contents.html}.

\bibitem{we:almost}
Weinstein, A., Almost invariant submanifolds for compact group
actions, preprint math.DG/9908133.

\bibitem{we:proper}
Weinstein, A., Proper groupoids (in preparation).

\end{thebibliography}
\end{document}